\documentstyle[11pt,righttag]{amsart}

\begin{document}

\topmargin-0.1in
\textheight8.5in
\textwidth5.5in
 
\footskip35pt
\oddsidemargin.5in
\evensidemargin.5in

\newcommand{\V}{{\cal V}}      % cal in math mode
\renewcommand{\O}{{\cal O}}
\newcommand{\LL}{\cal L}
\newcommand{\Ext}{\hbox{Ext}}
\newcommand{\pE}{\hbox{$^\pi\kern-2pt E$}}
\newcommand{\hQ}{\hbox{$\hat Q$}}
\newcommand{\phQ}{\hbox{$ '{\hat Q}$}}
\newcommand{\phd}{\hbox{$ '{\hat \delta}$}}

\newcommand{\lonto}{{\protect \longrightarrow\!\!\!\!\!\!\!\!\longrightarrow}}

\newcommand{\m}{{\frak m}}
\newcommand{\gl}{{\frak g}{\frak l}}
\newcommand{\ssl}{{\frak s}{\frak l}}

\renewcommand{\d}{\partial}

\newcommand{\ds}{\displaystyle}
\newcommand{\s}{\sigma}
\renewcommand{\l}{\lambda}
\renewcommand{\a}{\alpha}
\renewcommand{\b}{\beta}
\newcommand{\G}{\Gamma}
\newcommand{\g}{\gamma}

\newcommand{\C}{{\Bbb C}}
\newcommand{\N}{{\Bbb N}}
\newcommand{\Z}{{\Bbb Z}}
\newcommand{\ZZ}{{\Bbb Z}}
\newcommand{\K}{{\cal K}}
\newcommand{\ve}{{\varepsilon}}
\newcommand{\cupp}{{\star}}

\newcommand{\rowxy}{(x\ y)}
\newcommand{\colxy}{ \left({\begin{array}{c} x \\ y \end{array}}\right)}
\newcommand{\scolxy}{\left({\begin{smallmatrix} x \\ y
\end{smallmatrix}}\right)}

\renewcommand{\P}{{\Bbb P}}

\newcommand{\la}{\langle}
\newcommand{\ra}{\rangle}

\newtheorem{thm}{Theorem}[section]
\newtheorem{lemma}[thm]{Lemma}
\newtheorem{cor}[thm]{Corollary}
\newtheorem{prop}[thm]{Proposition}

\theoremstyle{definition}
\newtheorem{defn}[thm]{Definition}
\newtheorem{notn}[thm]{Notation}
\newtheorem{ex}[thm]{Example}
\newtheorem{rmk}[thm]{Remark}
\newtheorem{rmks}[thm]{Remarks}
\newtheorem{note}[thm]{Note}
\newtheorem{example}[thm]{Example}
\newtheorem{problem}[thm]{Problem}
\newtheorem{ques}[thm]{Question}
 
\numberwithin{equation}{section}%this numbers equations by section

\renewcommand{\proofname}{\bf \rom{Proof}}
\newcommand{\onto}{{\protect \rightarrow\!\!\!\!\!\rightarrow}}
\newcommand{\donto}{\put(0,-2){$|$}\put(-1.3,-12){$\downarrow$}{\put(-1.3,-14.5) 

{$\downarrow$}}}

\newcounter{letter}
\renewcommand{\theletter}{\rom{(}\alph{letter}\rom{)}}

\newenvironment{lcase}{\begin{list}{~~~~\theletter} {\usecounter{letter}
\setlength{\labelwidth4ex}{\leftmargin6ex}}}{\end{list}}

\newcounter{rnum}
\renewcommand{\thernum}{\rom{(}\roman{rnum}\rom{)}}

\newenvironment{lnum}{\begin{list}{~~~~\thernum}{\usecounter{rnum}
\setlength{\labelwidth4ex}{\leftmargin6ex}}}{\end{list}}

%\begin{document}
 
\title{Generalizing the notion of  Koszul Algebra}


\author[Cassidy and Shelton]{ }
\maketitle

\begin{center}

\vskip-.2in Thomas Cassidy$^\dagger$ and Brad Shelton$^\ddagger$\\
\bigskip

$^\dagger$Department of Mathematics\\ Bucknell University\\
Lewisburg, Pennsylvania  17837
\\ \ \\

   $^\ddagger$Department of Mathematics\\ University of Oregon\\
Eugene, Oregon  97403-1222
\end{center}

\setcounter{page}{1}

\thispagestyle{empty}

 \vspace{0.2in}

\begin{abstract}
\baselineskip15pt
We introduce a generalization of the notion of a Koszul algebra, which includes graded algebras with relations in different degrees, and we establish some of the basic properties of these algebras.  This class is closed under twists, twisted tensor products, regular central extensions and Ore extensions. We explore the monomial algebras in this class and we include some well-known examples of algebras that fall into this class.   \end{abstract}

\bigskip

\baselineskip18pt

%%%%%%%%%%%%%%%%%%%%%%%%%%%%%%%%%%%%%%%%%%%%%%%%%%%%%%%%%%%%%%%%%%%%%%
\section{Introduction}
Koszul algebras were originally defined by Priddy in 1970 \cite{Priddy}  and have since revealed important applications in algebraic geometry, Lie theory, quantum groups,  algebraic topology and, recently, combinatorics (cf. \cite{LH}).    The rich structure and long history of Koszul algebras  are clearly detailed in \cite{PP}.

There exist numerous equivalent definitions of a Koszul algebra  (see for example \cite{BackFro}).  For our purposes,  there are two particular equivalent definitions that motivate our discussion. Let $K$ be a field and let $A$ be connected graded $K$-algebra, finitely generated in degree 1.  Let $E(A) = \bigoplus E^{n,m}(A) = \bigoplus Ext_A^{n,m}(K,K)$ be the associated bigraded Yoneda algebra of $A$ (where $n$ is the cohomology degree and $-m$ is the internal degree inherited from the grading on $A$).  Set $E^n(A) =\bigoplus_mE^{n,m}(A)$.  Then $A$ is said to be  Koszul  if it satisfies either of the following (equivalent) definitions:

(1)  {\it Diagonal purity}:  $E^{n,m}(A) = 0$ unless $n=m$.

(2) {\it Low-degree generation}:  $E(A)$ is generated as an algebra by $E^{1,1}(A)$. 

\noindent  Stated this way, these are very strong conditions and either of them immediately imply that a Koszul algebra must be quadratic. 
 
Recently, Berger introduced the class of $N$-Koszul algebras \cite{Berg}  as a possible generalization of Koszul algebras.  The definition of $N$-Koszul generalizes the purity definition of Koszul given above and is given as the statement:  $E^{n,m}(A)=0$ unless $m = \delta(n)$, where $\delta(n)=N(n-1)/2+1$ if $n$ is odd and $\delta(n) = Nn/2$ if $n$ is even.  In particular, an $N$-Koszul algebra must have all of it relations in degree $N$ and each $E^n(A)$ is {\it pure} in the sense that it is supported in a singe internal degree.  Purity then becomes a powerful homological tool, see for example \cite{BGZ} or \cite{GMMZ}.

In \cite{GMMZ} it was shown that $N$-Koszul could also be rephrased in terms of low-degree generation.  It is shown there, as we will reprove in \ref{N-Kos}, that $A$ is $N$-Koszul if and only if:
\begin{enumerate}
\item $E(A)$ is generated as an algebra by $E^1(A)$ and $E^2(A)$, and
\item $A$ has defining relations all of degree $N$.
\end{enumerate} 

Several works (\cite{Berg}, \cite{GMMZ}, \cite{BGZ}, \cite{FV} among others) demonstrate similarities between $N$-Koszul algebras and Koszul algebras, and so provide evidence that the $N$-Koszul algebras should be included in a generalization of Koszulity.  However the restriction to $N$-homogeneity is somewhat artificial and poses serious problems.  In particular, the class of $N$-Koszul algebras, for $N>2$, is not closed with respect to many standard operations, such as graded Ore extensions, regular normal extensions, or tensor products. 

The goal of this paper is to give an alternate generalization of the notion of Koszul with four good properties: it should allow for graded algebras with relations in more than one degree; it should, at minimum, be closed under tensor products and regular central extensions; it should include the $N$-Koszul algebras, and it should collapse to the definition of Koszul in the case of quadratic algebras.  As we will see, the following simple definition does all this, and more. 

\begin{defn}  The graded algebra $A$ is said to be $\K_2$ if $E(A)$ is generated as an algebra by 
$E^1(A)$ and $E^2(A)$.
\end{defn}  
It is clear that this is the next most restrictive definition one could make, following Koszul and $N$-Koszul, since for a non-Koszul algebra, $E(A)$ could never be generated by anything less than $E^1(A)$ and $E^2(A)$.  However, this definition sacrifices homological purity.  Surprisingly, many statements that one would want in a generalization of Koszulity hold even without the convenience of homological purity. 

The same definition, under the name {\it semi-Koszul}, was used  in \cite{Mauger} to study certain Hopf algebras over finite fields.   Clearly this definition captures the Koszul and $N$-Koszul algebras.   But there are many important examples of $\K_2$ algebras which have defining relations in more than one degree and so cannot be $N$-Koszul.   For example, the homogeneous coordinate ring of any projective complete intersection is $\K_2$ (Corollary \ref{complete}).   Any Artin-Schelter regular algebra of global dimension four on three generators (refer to \cite{ATV1} or \cite{LPWZ}) will have two quadratic relations and two cubic relations.  We show in Theorem \ref{AS} that these algebras are all $\K_2$.  

The main results in this paper describe operations that will preserve the class of $\K_2$ algebras.  For example, the class is closed under twists and twisted tensor products, as proved in section 3.  Most of the results require technical hypotheses, and are therefore not stated explicitly in this introduction.  Along the way we give an alternative description of $\K_2$ algebras in terms of minimal projective resolutions (Theorem \ref{MC}) and define $\K_2$-modules.  This is analogous to the standard Koszul theory of modules with linear projective resolutions.   As a corollary to this result we get Theorem \ref{monomialthm}, a simple combinatorial algorithm for establishing which monomial algebras are $\K_2$.

The last five sections of the paper deal with change of rings theorems.  Section 6 contains lemmas about the structure of the spectral sequences associated to a change of rings.  In Section 7 we use the notion of a $\K_2$ module and prove Theorem 
\ref{ringchange} about algebras with $\K_2$ quotients, an analog to a standard theorem in Koszul theory. In particular, we see then that the class of $\K_2$ algebras is closed under regular central extensions.  In Theorem \ref{quotient>1} we analyze the inheritance of the $\K_2$ property from an algebra $A$ to a normal regular quotient of the form $A/gA$ where the degree of $g$ is 1.    The case where $g$ has degree greater than 1 is analyzed in Theorem \ref{normalreg}, with a remarkably different answer.  In Theorem \ref{Ore} we show that the class of $\K_2$ algebras is closed under graded Ore extension ($A=B[z,\sigma;\delta]$), and in the process we establish an interesting and apparently previously unknown fact about maps induced on the Yoneda algebra (Lemma
\ref{cool}). 

Throughout the paper we include many examples of both $\K_2$ and non-$K_2$ algebras that we hope are both interesting and illustrative.

%%%%%%%%%%%%%%%%%%%%%%%%%%%%%%%%%%%%%%%%%%%%%%%%%%%%%%%%%%%%%%%%%%%%%%
\section{Notation}

Let $K$ be a field, $V$ a finite dimensional $K$-vector space and $T(V)$ the usual $\N$-graded tensor algebra over $V$.  Throughout this paper, the phrase {\it graded algebra} will refer to a connected graded $K$-algebras of the form $A=T(V)/I$, where
$I$ is an ideal of $T(V)$ that is generated by a (minimal) finite collection of homogeneous elements 
$\{r_1\ldots, r_s\}$ of degree at least 2.  We do not asssume that the relations $r_i$ all have the same degree. 

   We put
$I' = I\otimes V+V\otimes I$ so that $I'_n$ is the set of element of $I$ of degree $n$ that are generated by elements of $I$ of strictly smaller degree.  We note that $I'_m = I_m $ for all $m$ sufficiently large, and that $I'$ is a proper sub-ideal of $I$.

\begin{defn}\label{essential}
An element $r\in I$ is said to be {\it essential} (in $I$) if $r$ is not in $I'$, or equivalently $r$ does not vanish in $I/I'$.
\end{defn}  

For any $\ZZ$-graded vector space $W$, we use $W(n)$ to denote the same vector space with degree shifted  grading $W(n)_k = W_{k+n}$.   It is sometimes convenient to write $W(n_1,n_2,\ldots,n_k)$ for
$\bigoplus_{j=1}^{j=k} W(n_j)$.   We denote the projective dimension of $W$ as an $A$-module by $pd_A(W)$.

Given a graded algebra $A$, as above, we identify the ground field $K$ with the trivial (left or right) $A$-module $A/A_+$ and we write $\epsilon: A \to K$ for the corresponding augmentation. 

Our purpose is to study certain properties of the bigraded Yoneda algebra of $A$ which we denote
$E(A)$.  This is defined as $E(A) = \bigoplus_{n,k\ge 0} E^{n,k}(A)$ where 
$E^{n,k}(A) = Ext_A^{n,k}(K,K) = Ext_A^n(K,K)_{-k}$.  In the bigrading $E^{n,k}(A)$, $n$ is the {\it cohomology} degree and $-k$ is the {\it internal} degree (inherited from the grading on $A$).   Note that $E(A)$ is supported in non-positive internal degrees.  We also write $E^n(A)$ for $\bigoplus_k E^{n,k}(A)$.  The cup product on the ring $E(A)$ is denoted
$\cupp$.   The $Ext$ groups involved here are calculated in the category of  $\Z$-graded
and locally finite $A$-modules with graded $Hom$ spaces.  

Given two graded algebras $A$ and $B$ and a degree zero $K$-algebra homomorphism 
$\sigma:A \to B$, we will denote the induced algebra homorphism from $E(B)$ to $E(A)$ by 
$\sigma^*:E(B) \to E(A)$.  In particular, if $\sigma$ is a graded automorphism of $A$, then 
$\sigma^*$ is a bi-graded automorphism of $E(A)$.

%%%%%%%%%%%%%%%%%%%%%%%%%%%%%%%%%%%%%%%%%%%%%%%%%%%%%%%%%%%%%%%%%%%%%%
\section{Operations on the class $\K_2$}

In this section we show that the class of $\K_2$ algebras is closed under twists by automorphisms and twisted tensor products. 

Let $\sigma$ be an automorphism of a graded algebra $A = T(V)/I$. Then we have the notion of the twisted algebra $A^\sigma$, as studied extensively by 
\cite{Zhang-twist}.  This algebra has the same underlying graded vector space structure as $A$, with the product of homogeneous elements $a\cdot b = a\sigma^{|a|}(b)$.  (Here
$ |a|$ gives the degree of $a$.)  

\begin{thm}\label{twist} Let $A$ be an $\N$-graded algebra and $\sigma$ a graded automorphism of $A$.  Then $A$ is $\K_2$ if and only if $A^\sigma$ is $\K_2$.
\end{thm}

\begin{pf}
The Theorem follows immediately from the claim  that 
$E(A^\sigma) = E(A)^{(\sigma^{-1})^*}$, where the twisting is done with respect to the internal grading.  We sketch a proof of this claim.   The spaces in the bar complex computing 
$E(A)$ and $E(A^\sigma)$ are the same, namely 
$\Omega_n:=Hom_K(A_+^{\otimes n},K)$.  Let $d$ be the differential on $\Omega_n$ which determines $E(A)$, and 
$d_\sigma$ the differential which determines $E(A^\sigma)$.  Then there is an explicit homotopy  $\mu^*:(\Omega_n,d)\to (\Omega_n,d_\sigma)$ between the two complexes.  This homotopy is the dual of a vector space isomorphism $\mu\in Aut_K(A_+^{\otimes n})$ with the following complicated formula.  Let $\bar a = a_1\otimes \ldots \otimes a_n$ be a multihomogeneous element of $(A_+)^{\otimes n}$.  Define $|\bar a,i| = \sum_{j=i}^n |a_j|$ for $1\le i\le n$.  Then we define:
$$\mu(a_1\otimes \ldots\otimes a_n)= \mu(\bar a) = 
	\sigma^{|\bar a,1|}(a_1)\otimes \sigma^{|\bar a,2|}(a_2) \ldots 
			\otimes \sigma^{|\bar a,n|}(a_n).$$
One checks that $d(f\circ \mu) =(d_\sigma f)\circ \mu$.   This shows that $\mu^*$ induces a bigraded vector space isomorphism from $E(A)$ to $E(A^\sigma)$.  Let $\cupp$ be the usual cup product in $E(A)$ and $\cupp_\sigma$ the cup product in $E(A^\sigma)$. One then calculates explicitly that for bigraded elements $f$ and $g$ of $E(A)$
$$\mu^*f \cupp_\sigma \mu^* g = \mu^*( (\sigma^*)^{-|g|}(f) \cupp g )$$
where $|g|$ is the internal degree of $g$. 
This shows that $E(A^\sigma)$ is isomorphic to $E(A)^{(\sigma^{-1})^*}$ and proves the claim.
\end{pf}

Let $A$ and $B$ be graded algebras and let $\sigma$ be a graded automorphism of $A$.  Then we can define an associative multiplication on the $K$-tensor product $B\otimes A$ by the rule $(b_1\otimes a_1)\cdot (b_2\otimes a_2) = b_1b_2 \otimes \sigma^{|b_2|}(a_1)a_2$ whenever $a_1, a_2, b_1, b_2 $ are homogeneous.  We denote the graded algebra thus formed by $B\otimes^\sigma A$ and write
 $b\otimes^\sigma a$ for the element $b\otimes a$ in the algebra.

\begin{thm}  If $A$ and $B$ are graded algebras,  and $\sigma$ is an automorphism of $A$, then $B\otimes^\sigma A$ is $\K_2$ if and only if $A$ and $B$ are both $\K_2$.
\end{thm}

\begin{pf}      Arguing in a similar fashion to the previous Theorem, one proves the formula:
$$E(B\otimes^\sigma A) = E(B) \overline\otimes^{(\sigma^{-1})^*} E(A).$$
(See also \cite{PP}, 3.1.1.)  The result follows immediately.  
\end{pf}

%%%%%%%%%%%%%%%%%%%%%%%%%%%%%%%%%%%%%%%%%%%%%%%%%%%%%%%%%%%%%%%%%%%%%%
\section{Conditions on Minimal Projective Resolutions}

It is well known that a graded algebra $A$ is Koszul if and only if the trivial module admits a linear free resolution, see for example \cite{PP}.  A simple way to express this is to say that all of the matrix entries of the maps in any minimal projective resolution of the trivial $A$-module $K$  have degree 1.  In this section we prove a similar, but much more complicated, version of this statement for $\K_2$ algebras, and then we use this to determine which 4 dimensional Artin-Schelter regular algebras are $\K_2$.  

We continue with the notation $A=T(V)/I$ where $I$ is a graded ideal generated (minimally) by homogeneous elements $r_1,r_2,\ldots,r_s$ of degrees $n_1,\ldots ,n_s$ respectively (all greater than 1).  Let $R$ denote both the set $\{r_1,\ldots,r_s\}$ and the transpose of the $1\times s$ matrix $(r_1\ r_2\ \ldots\  r_s)$. We fix a minimal projective resolution
$(Q,d)$ of the trivial left $A$ module $K$.  It is clear that if $A$ is a $\K_2$ algebra then
every $Q^n$ is finitely generated over $A$. Hence we will assume unless otherwise stated that every $Q^n$ is finitely generated.  For each $n\ge0$, we can define integers $t_n$ and $m_{1,n}\cdots,m_{t_n,n}$ by
$Q^n = A(m_{1,n}\cdots,m_{t_n,n})$.  (Note that $m_{k,n} \le -n$ for all $k$.)

We choose a homogeneous $A$-basis for each $Q^n$.  With respect to these bases, we let $M_n$ be the matrix representing $d_n:Q^n\to Q^{n-1}$ and let 
$\hat M_n$ be a lift of $M_n$ to a matrix over $T(V)$ with homogeneous entries.  
The entries of  $\hat M_1$ give a basis for $V$ and we may choose our basis of $Q^2$, as well as the lift $\hat M_2$, so that 
$\hat M_2\hat M_1 = R$.  (When multiplying matrices defined over $T(V)$ we suppress the tensor product notation.)

The following Lemma is obvious.

\begin{lemma}\label{obvious}  Let $\alpha_1,\ldots,\alpha_s$ be homogeneous elements of 
$T(V)$ and suppose $r = \sum_s \alpha_s r_s$.  Then $r$ is essential in $I$ if and only if some 
$\alpha_s$ is a unit. 
\end{lemma}

Before we state our theorem we recall the minimal projective resolution contruction of the cup product on $E(A)$.  Recall that $E^n(A) = Hom_A(Q^n,K)$, by minimality of the resolution. To multiply $f\in E^n(A)$ by $g\in E^k(A)$ we consider the following commutative diagram, where existence of the downward maps $F_{j}$ is assured by projectivity.

$$\begin{array}{ccccccclc}
Q^{n+k} &\kern -12 pt \buildrel M_{n+k} \over \to &\kern -16pt Q^{n+k-1} &\kern -16pt \to  \ldots\to  & Q^{n+1}
					&\kern -12pt \buildrel M_{n+1} \over \to &\kern -16 pt Q^n \\
 F_{k}\downarrow\hphantom{F_{k}} &&
	 F_{k-1}\downarrow\hphantom{F_{k-1}}&&   
	 F_{1}\downarrow\hphantom{F_{1}} &&
	 F_{0}\downarrow\hphantom{F_{0}}&
	\kern -12 pt\buildrel f\over \searrow\\[7 pt] 
Q^{k} &\kern -12 pt \buildrel M_{k} \over \to &\kern -16 pt Q^{k-1} &\kern-16pt  \to  \ldots\to  & Q^{1}
					&\kern -12pt \buildrel M_{1} \over \to &\kern -16 pt Q^0 &
						\kern -12 pt \buildrel \epsilon \over \to & K \\[7 pt]
g \downarrow \hphantom{g}\\[7pt]
K\\
					\end{array}
$$
Then the product of $g$ and $f$ in $E(A)$ is 
$g \cupp f := g\circ F_k \in Hom_A(Q^{n+k},K)$. 

It will be convenient to work with a graded basis of $E(A)$.  Let $\{\ve_{k,n}\}^{t_n}_{k=1}$ with $\ve_{k,n}:Q^n\to K$ be the graded basis of $E^n(A)$ dual to a minimal set of homogeneous $A$-generators of $Q^n$, so that $\ve_{k,n}$ has cohomology degree $n$ and internal degree $m_{k,n}$ (since the corresponding $A$-generator has degree 
$-m_{k,n}$).
 
 For $n\ge 1$  let
$U_n=E^2(A) \cupp E^{n-2}(A)$ and $V_n=E^1(A)\cupp E^{n-1}(A)$.  We note that the algebra $A$ is $\K_2$ if and only if $E^n(A) = U_n + V_n$ for all $n\ge 3$.  (This requires an inductive argument.)

We set some temporary notation to analyze the diagrams dictating the relationship of $E^n(A)$ to $U_n$ and to $V_n$. 

For every  $n$ and every $1\le k\le t_n$, let $e_{k,n}$ be the $k^{th}$ column of the 
$t_n\times t_n$ identity matrix.   We can define matrices $G_{k,n}$ and $J_{k,n}$  with entries in $A$ by commutativity of the following diagram.
$$\begin{array}{ccccclc}
Q^{n+2} &\kern -16pt \buildrel M_{n+2} \over \to &\kern -16pt Q^{n+1} &\kern -16pt \buildrel M_{n+1}\over\to &\kern -16pt Q^{n} \\
 J_{k,n}\downarrow\hphantom{J_{k,n}} &&\kern -16pt
	 G_{k,n}\downarrow\hphantom{E_{k,n}}&&\kern -16pt   
	 e_{k,n}\downarrow\hphantom{e_{k,n}} &\kern -16pt
	\buildrel \ve_{k,n} \over \searrow\\[7 pt] 
Q^{2} &\kern -16pt \buildrel M_{2} \over \to &\kern -16pt Q^{1} &\kern -16pt\buildrel M_1 \over \to 
					&\kern -16pt A
					& \kern -16pt\buildrel \epsilon \over \to & K \\
					\end{array}
					\hfill(*)
$$
The entries of $J_{k,n}$ and $G_{k,n}$ are all homogeneous, but the actual degrees will not be important.  Since the matrix $e_{k,n}$ is scalar and the matrix $M_1$ is linear, we
may lift the matrix $G_{k,n}$ to a (unique) homogeneous matrix $\hat G_{k,n}$ over 
$T(V)$ which satisfies $\hat G_{k,n}\hat M_1 = \hat M_{n+1}e_{k,n}$.  We also choose a homogeneous lift
of $J_{k,n}$ to $\hat J_{k,n}$ in $T(V)$.  Let $R_{k,n} = \hat M_{n+2}\hat G_{k,n} - \hat J_{k,n}\hat M_2$. We observe that $R_{k,n}$ is a matrix of homogeneous elements of $I$.  

With all of this notation in hand we can state two technical lemmas, the first of which concerns  $U_n$.

\begin{lemma}\label{cup1} (1) Let $1\le p\le t_2$ and $1\le k\le t_n$ and define scalars $\lambda_j$ by 
$\ve_{p,2}\cupp \ve_{k,n} = \sum_{j=1}^{t_{n+2}} \lambda_j \ve_{j,n+2}$.  Then
for $1\le i\le t_{n+2}$, $\lambda_i \ne 0$ if and only if the entry in row $i$ and column 
$p$ of $J_{k,n}$ is a unit.  

(2) Let $1\le i\le t_{n+2}$.  Then there exist $p$ and $k$ as in (1) such that $\lambda_i\ne 0$ if and only if some entry in row $i$ of $\hat M_{n+2} \hat M_{n+1}$ is essential in $I$.  
\end{lemma}

\begin{pf}

The first statement is clear (and the required unit is $\lambda_i$).  

For the second statement, suppose first that the required 
$p$ and $k$ exist.  Then the entry in row $i$ and column $p$ of $J_{k,n}$ is the unit
$\lambda_i$.  This must also be true of the matrix $\hat J_{k,n}$.  Let $f_i$ be row $i$ of the $t_{n+2}\times t_{n+2}$ identity matrix.  Then $f_i\hat J_{k,n}R := r$ is essential in $I$ by \ref{obvious}.
But we have 
$$r = f_i(\hat M_{n+2}\hat G_{k,n}-R_{k,n})\hat M_1 = f_i\hat M_{n+2}\hat M_{n+1} e_{k,n} - f_i R_{k,n}\hat M_1.$$
Since $R_{k,n}$ has entries from $I$, the element $f_iR_{k,n}\hat M_1$ cannot be essential.  It follows that 
$f_i\hat M_{n+2}\hat M_{n+1}e_{k,n}$ is essential, as required.

Conversely, suppose that an element of row $i$ of $\hat M_{n+2}\hat M_{n+1}$ is essential, say the element in column $k$. Call the element $r$.  Then 
$$r=f_i\hat M_{n+2}\hat M_{n+1}e_{k,n} = f_i(\hat J_{k,n}\hat M_2  + R_{k,n})\hat M_1 .$$
Since no entry of $f_iR_{k,n}\hat M_1$ is essential, we see that
$ f_i\hat J_{k,n}\hat M_2\hat M_1 = f_iJ_{k,n}R$ is essential.  From Lemma \ref{obvious}, some entry of $f_iJ_{k,n}$ is thus a unit, as required.
\end{pf}

Our second technical lemma concerns $V_n$ and we leave its (easier) proof to the reader.

\begin{lemma}\label{cup2} (1) Let $1\le p\le t_1$ and $1\le k\le t_n$ and define scalars $\mu_j$ by
$\ve_{p,1}\cupp \ve_{k,n} = \sum_{j=1}^{t_{n+1}} \mu_j \ve_{j,n+1}$.  Then
for $1\le i\le t_{n+1}$, $\mu_i \ne 0$ if and only if the entry in row $i$ and column 
$p$ of $E_{k,n}$ is a unit.  

(2) Let $1\le i\le t_{n+1}$.  Then there exist $p$ and $k$ as in (1) such that $\mu_i\ne 0$ if and only if some entry in row $i$ of $\hat M_{n+1}$ is (nonzero) linear. 

\end{lemma}

To state our theorem we need just a bit more notation.  For each $n\ge 2$, 
let $L_n$ be the image of $\hat M_n$ modulo the ideal $T(V)_{\ge 2}$ and
let $E_n$ be the image of $\hat M_n\hat M_{n-1}$ modulo $I'$.  We think of $L_n$ as the linear part of $\hat M_n$ and $E_n$ as the essential part of $\hat M_n\hat M_{n-1}$.
Finally, let $[L_n:E_n]$ be the $t_n\times (t_{n-1}+t_{n-2})$ matrix obtained by concatenating the rows of $L_n$ and $E_n$.  
(Note that the entries in any given column of this matrix are all in the same vector space, either $V$ or $I/I'$.)

\begin{thm}\label{MC}  Let $A$ be a graded $K$-algebra as above, and $Q\to K \to 0$ a minimal graded projective resolution.  Then the following are equivalent:

(1) $A$ is $\K_2$ 

(2) For  $2<n\le pd_A(K)$, $Q^n$ is finitely generated and the rows of the matrix $[L_n:E_n]$ are linearly independent over $K$.  
\end{thm}

\begin{pf}

Suppose first that for some specific $n>2$, the rows of the matrix 
$[L_n:E_n]$ are linearly dependent.   By changing basis (homogeneously) in the free module $Q^n = A(m_{1,n},\ldots, m_{t_n,n})$, we may assume that the first row of the matrix $[L_n:E_n]$ is zero. By lemmas \ref{cup1} and \ref{cup2}, the coefficient of the cohomology class $\ve_{i,n}$ in every element of $U_n$ or $V_n$ is  then zero.  Hence $\ve_{i,n} \not\in U_n + V_n$ and $A$ is not $\K_2$.

Conversely, suppose that $A$ is not $\K_2$ and fix the smallest $n$ for 
which $E^n(A) \ne U_n + V_n$.  Choose a homogeneous basis for 
$U_n + V_n$, say
$\ve_{1,n},\ldots,\ve_{k,n}$, $k<t_n$.  Extend this a to full homogeneous basis of $E^n(A)$ by choosing
$\ve_{k+1,n},\ldots,\ve_{t_n,n}$.  This basis then corresponds to a homogeneous $A$-basis for $Q^n$.  We may assume the matrices representing the minimal resolution are calculated with respect to this new basis.  But then again by lemmas 
\ref{cup1} and \ref{cup2},  the last row of $[L_n:E_n]$ must be zero.  
\end{pf} 

 An algebra will often fail the $\K_2$ hypothesis simply because of the existence of a cohomology class in the ``wrong'' internal degree, i.e. an internal degree that could not be generated by the internal degrees of elements of lower cohomological degree.  
For example, the algebra
$B=K\langle x,y\rangle / \langle x^2-xy, y^2\rangle$ is a well-known  example of a non-Koszul quadratic (and hence non-$\K_2$ algebra).  Since the algebra is quadratic, the portion of $E(A)$ generated by $E^1(A)$ and $E^2(A)$ is exactly $\sum_nE^{n,n}(A)$.   But $E^{3,4}(B) \ne 0$, so we have cohomology in the ``wrong'' internal degree.  

However, an algebra can also fail to be $\K_2$ even when the errant cohomology class has an internal degree that could be generated by classes of lower cohomological degree.  Theorem \ref{MC} can be 
a simple way to detect this, as exhibited by the following example.

\begin{example}   Let $A$ be the algebra $K\langle x,y\rangle / \langle x^2-xy, yx, y^3\rangle$.  Using the fact that $A$ has Hilbert series $1+2t+2t^2$, one can calculate the first few terms of a minimal projective resolution of $_AK$ to obtain:

$$A(-3,-4,-4,-4,-4) \buildrel M_3\over \to A(-2,-2,-3) \buildrel M_2\over \to A(-1,-1) \buildrel M_1\over \to A\to K\to 0$$
where
$$
\hat M_3= \left(\begin{array}{ccc}y&0&0\\ 0&xy&0\\ 0&y^2&0\\ 0&0&x\\ 0&0&y  \end{array}
\right) \qquad 
\hat M_2 = \left(\begin{array}{cc} x & -x\\ y & 0 \\  0& y^2 \end{array}\right) \qquad
\hat M_1= \left(\begin{array}{c} x \\ y \end{array}\right) 
$$
Since the tensor $xy^2$ is in the ideal of definition of $A$, but is not essential ($xy^2 = -x(x^2-xy) + (x^2-xy)(x-y) + x(yx)$),
we get
$$[L_3:E_3] = \left[\begin{array}{ccccc}y&0&0&yx&-yx \\ 0&0&0&0&0\\ 0&0&0&y^3&0\\ 0&0&x&0&0\\
0&0&y&0 & y^3 \\ \end{array}\right]$$
By \ref{MC}, the algebra $A$ is not $\K_2$ and  $E^3(A) \ne U_3 +V_3$.  Indeed a more careful analysis of the  matrix above shows $\dim_K((U_3+V_3)_{-4}) = 3$, whereas we know
$\dim_K(E^{3,4}(A)) = 4$.  

\end{example}

We can also reprove the following fact about $N$-Koszul agebras (see also \cite{GMMZ}). 

\begin{cor}\label{N-Kos}  The algebra $A$ is $N$-Koszul if and only if it is $N$-homogeneous and $\K_2$.
\end{cor}

\begin{pf}
Suppose first that the algebra $A$ is $N$-Koszul.  Then the matrices $\hat M_n$ will have linear entries when $n$ is odd and entries all of degree $N-1$ when $n$ is even.  In particular $\hat M_n\hat M_{n-1}$ will always be a matrix whose entries are all essential relations.  The theorem immediately implies that $A$ is $\K_2$.

Conversely, suppose that $A$ is $\K_2$ and $N$-homogeneous.  Then a relation of $A$ is essential if and only if it has degree $N$.  An inductive argument, based on the condition of  Theorem \ref{MC}, then assures that the matrices $\hat M_n$ must alternate between having entries of degree 1 when $n$ is odd and degree $N-1$ when $n$ is even.  This proves that $A$ is $N$-Koszul.
\end{pf}

We conclude this section by determining which 4 dimensional Artin-Schelter regular (AS-regular) algebras are $\K_2$.  See \cite{LPWZ} for the definition of Artin-Schelter regular and for details on the degrees of the defining relations for these algebras.  
 
Let $A$ be an AS-regular algebra of global dimension 4.  Then $A$ has 2, 3, or 4 linear generators.   If $A$ has 4 generators then it is Koszul \cite{ST} and hence $\K_2$.  If the algebra has 2 generators then it has two relations, one of degree 3 and one of degree 4.  The definition of Artin-Schelter regular assures us that $_AK$ has a minimal projective resolution of the form:$$
0 \to A(-7)   \to A(-6)^2    \to A(-3,-4) 
 \to A(-1)^2     \to A \to K \to 0
$$
and it is easy to see that the algebra can not be $\K_2$ as there is a nonzero cohomology class in 
$E^{3,6}(A)$ which can not be generated by $E^2(A) = E^{2,3}(A) + E^{2,4}(A)$ and 
$E^1(A) = E^{1,1}(A)$.  

For the remainder of this section we will assume that $A$ is AS-regular of global dimension 4 and has 3 linear generators.  Such an algebra will have 4 relations of degrees 2, 2, 3 and 3.   We will write $A$ as $K\langle x_1,x_2,x_3\rangle/I$.   In this case $_AK$ has a minimal projective resolution of the form:
\begin{equation}\label{ASresolution}
0 \to A(-5) \buildrel Y \over \to A(-4)^3 \buildrel N \over \to A(-2,-2,-3,-3) 
\buildrel M \over \to A(-1)^3 \buildrel X \over \to A \to K \to 0
\end{equation}
where the matrices $Y$, $N$, $M$, $X$ have homogeneous entries from 
$K\langle x_1,x_2,x_3\rangle$.  From the Gorenstein condition, we may assume:
$$
Y =\left(\begin{array}{ccc}x_1 &x_2 & x_3\end{array}\right),\quad
X=\left(\begin{array}{c}x_1 \\x_2 \\ x_3\end{array}\right),
$$
$N$ is a 3 by 4 matrix where the entries of the first two columns are 2-tensors and those in the last two columns are 1-tensors and $M$ is a 4 by 3 matrix with 1-tensors as entries in the first two rows and 2-tensors in the last two rows.  In particular, the entries in $MX$ and the entries in $YN$ must each give a set of relations for $A$.  The $K$-span of the two quadratic entries in $MX$ must be
the same as the $K$-span of the two quadratic entries in $YN$.

\begin{thm}\label{AS}  Let $A$ be an AS-regular algebra of global dimension 4 with 3 linear generators.   Then $A$ is $\K_2$.
\end{thm}

\begin{pf}
We see from \eqref{ASresolution} that the $\K_2$ matrix condition given in Theorem \ref{MC}  can fail only at $N$.  Let us assume that the condition fails, and find a contradiction.

Let $L$ be the linear part of the matrix $N$ and let $E$ be $NM$ modulo $I'$.  Failure of the matrix condition is the statement that the composite matrix
$[L:E]$ has linearly dependent rows.  By choosing new coordinates, we may assume that the first row of this matrix is zero.  Using this, let
$(\sum_i x_i \otimes a_i, \sum_i x_i\otimes b_i , 0, 0)$ be the first row of $N$, where 
the terms $a_i$ and $b_i$ are linear.

The two quadratic entries in $YN$ must span $I_2$, so from the first row of $N$ we see that no element of $I_2$ can contain a term beginning with $x_1$, i.e. if $\sum_i x_i \otimes l_i \in I_2$, then $l_1=0$.    
 
Consequently,  if $\sum_i x_i \otimes q_i \in (I')_3$, then 
$x_1\otimes q_1 \in (I')_3$ and $q_1 \in I_2$ (one can make no such statement about 
$q_2$ or $q_3$).  

Let the first two rows of $M$ form the matrix $(d_{i j})$, $i=1,2$ and $j=1,2,3$ where the elements $d_{i j}$ are all linear. Then the $j$-th entry of the first row of the matrix 
$NM$ is
$\sum_i x_i \otimes (a_i\otimes d_{1j}+ b_i \otimes d_{2j})$.   The statement that the first row of $E$ is zero is the statement that these three 3-tensors are all in $(I')_3$.  As 
stated above, this implies that 
$a_1\otimes d_{1,j} + b_1\otimes d_{2,j}$ is in $I_2$ for $1\le j\le 3$, so that $(a_1\ b_1\ 0\ 0)$ is in the kernel of $M$.  Since \eqref{ASresolution} is exact,  $(a_1\ b_1\ 0\ 0)$ must be in the image of $N$, which is impossible because the first two columns of $N$ are made up of two-tensors,  whereas $a_1$ and $b_1$ are linear.
\end{pf}

%%%%%%%%%%%%%%%%%%%%%%%%%%%%%%%%%%%%%%%%%%%%%%%%%%%%%%%%%%%%%%%%%%%%%%
\section{Monomial $\K_2$ Algebras}\label{monom}
  
In this section we explore in detail the $\K_2$ property for monomial algebras and using Theorem \ref{MC} we present an algorithm for determining whether a given monomial algebra is $\K_2$.

Let  $\{x_1,...,x_n\}$ be a fixed basis of the vector space $V$.  We use this basis to
identify $T=T(V)$ with $K\langle x_1,...,x_n\rangle$.  Throughout this section, $A = T/I$, where $I$ is an ideal generated by a minimal finite set of monic monomials (in the variables $x_i$) $R=\{r_1,..,r_s\}$.  For any element $a$ in $T$, we denote its image in 
$A$ by $\bar a$.   Let $T_{mon}$ be the set of monic monomials in $T$.  We write $|a|$ for the degree of any $a\in T_{mon}$. 

The combinatorics of monomial algebras is straightforward and built primarily on the finite set 
$LE(R)$ consisting of all $a\in T_{mon}$ for which there exists a $b\in T_{mon}$ such that $|b|>0$  and  $ab\in R$.
This is the set of proper ``left ends'' of $R$. 

For any pair of monomials $a$ and $b$ in $T_{mon}$ such that $\overline{ab}=0$ we define the monic monomial $L(a,b)$  to be $a'$ where $a=a''a'$ and $a'$ is minimal such that  $\overline{a'b}=0$.  If $\overline{ab}\ne 0$ we set $L(a,b)=0$.  Suppose that $\bar a\ne 0$ and $\bar b\ne 0$, but $\overline{ab}=0$.  Let $L(a,b)=a'$. Then there exists a minimal $b'$ so that $b=b'b''$ and $\overline{a'b'}=0$. The minimality of both $a'$ and $b'$ assures us that $a'b'\in R$, and hence that $a'\in LE(R)$.  

Finally, for $b\in T_{mon}$ such that $\bar b\ne0$, define 
 $$L(b)= \{a\in LE(R)\ |\  L(a,b)=a\}.$$
The following Lemma is clear from the definitions of $L(a,b)$, $L(b)$ and $LE(R)$.

\begin{lemma}\label{mon1} If $b$ is a monomial in $T$ and $\bar b\ne 0$ in $A$ then the left annihilator of $\bar b$ in $A$ is generated by $\{\bar a  |\ a\in L(b)\}$.  

\end{lemma}

From this Lemma it is easy to give a combinatorial description of the graded vector space structure of $E(A)$.  We do this by describing a minimal (monomial) projective resolution of the trivial module $_AK$ of the form:
$$ \to A^{t_m} \buildrel M_m\over \to A^{t_{m-1}} \to \cdots \to A^{t_1}\buildrel M_1\over \to A \to K\to 0.$$
To do this, it suffices to give an inductive definition of the matrices $\hat M_m$, the homogeneous lifts of the matrices $M_m$ to matrices over $T$. 

Let $\hat M_1 = \left(\begin{array}{c} x_1\\ \vdots \\ x_n\end{array}\right)$.  Assume that $m>1$ and that $\hat M_{m-1}$ has been defined and has the property that each row of 
$\hat M_{m-1}$ has a unique nonzero entry  and that entry is in $LE(R) \cup \{x_1,\ldots, x_n\}$.  Let
$\hat M_{m-1}$ have $t_{m-1}$ rows.  For $1\le i\le t_{m-1}$, if $b_i$ is the non-zero entry in row $i$ of $\hat M_{m-1}$ and 
$L(b_i) = \{a_1,\ldots,a_j\}$, let
$\hat M_m(i)$ be the $j$ by $t_{m-1}$ matrix with $a_1,\ldots, a_j$ arranged as column 
$i$ and zeroes elsewhere. If $L(b_i)$ is empty, this is the empty matrix. Finally, take 
$\hat M_m $ to be the matrix obtained by stacking the matrices $\hat M_m(i)$, $1\le i\le t_{m-1}$.  It is clear from construction that 
$\hat M_m$ has the property that each row has a unique nonzero entry from $LE(R) \cup \{x_1, \ldots, x_n\}$ (and in fact, for $m\ge 2$, all the nonzero entries are from $LE(R)$.)

We note that it is clear that $M_mM_{m-1} = 0$ for all $m>1$ and hence the matrices $M_m$ define a complex of free-graded $A$-modules as above.  Lemma \ref{mon1} makes the following clear.

\begin{lemma}  The complex described above is a minimal resolution of $_AK$.  

\end{lemma}

Now let $S$ be the smallest set in $T$ such that 
\begin{enumerate}
\item each $x_i$ is in $S$, and 
\item if $b\in S$ and $L(a,b)=a$ for some $a\in LE(R)$,  then $a\in S$.
\end{enumerate}
There is a simple inductive algorithm for calculating $S$.  Put $S_1=\{x_1\ldots,x_n\}$  and for $i>1$ define $S_i=\displaystyle \bigcup_{b\in S_{i-1}}L(b)\backslash\bigcup_{j<i}S_j$.   Then 
$S$ is the disjoint union of the $S_i$.  It should be clear that $S_k$ is just those nonzero matrix entries of $\hat M_k$ which have not appeared in $\hat M_i$ for any $i<k$.  Moreover, since $LE(R)$ is finite, $S_k = \emptyset $ for all $k$ sufficiently large.

\begin{thm}\label{monomialthm}
Let $A$ be a monomial algebra over the field $k$.  Then
$A$ is $\K_2$ if and only if the following condition holds for every $b\in S$:  if there exists an $a\in LE(R)$ such that $L(a,b)=a$ then $ab\in R$ or $deg(a)=1$. 
\end{thm}

\begin{pf}
From \ref{MC}, we need only see that the condition given in the theorem is equivalent to the condition that for all $m$
the rows of the matrix $[L_m:E_m]$ are linearly independent.  It is clear that the condition is equivalent to the statement that every row of $[L_m:E_m]$ is non-zero.  It therefore suffices to prove that 
$[L_m:E_m]$ has linearly dependent rows if and only if one of its rows is zero.  We begin with some simple observations about this matrix.

First we note that a monic {\it monomial} $r$ in $I$ is essential in $I$ if and only if $r\in R$, so we may assume 
that $E_m$ is a matrix with entries in the span of $R$ (rather than in $I/I'$).  Moreover, since every row of $\hat M_{m}\hat M_{m-1}$ has a single nonzero entry,  the nonzero rows of $E_m$ must be identical to the nonzero rows of $\hat M_{m}\hat M_{m-1}$.  It is now clear that if a set of rows of $E_m$ is linearly dependent and none of the rows are zero, then two of the rows must be identical.  We may assume that these are the first two rows and that the unique nonzero entry $r\in R$ is in the first column.  Now the multiplication
$\hat M_n \hat M_{n-1}$ must have the following form:
$$\left(\begin{array}{cccc}a & 0 & 0 &\cdots\\
				0 & b & 0 & \cdots \\
				\vdots
				\end{array}\right)
	\left(\begin{array}{cccc}c & 0 & 0 & \cdots\\
					    d & 0 & 0 & \cdots\\
					    \vdots 
					    \end{array}\right)
					    =
	\left(\begin{array}{cccc}r & 0 & \cdots\\
					    r & 0 & \cdots\\
					    \vdots
					    \end{array}\right)$$
where $a$, $b$, $c$ and $d$ are monomials with $ac=bd = r$. In particular there is a monomial $e$ such that $c = ed$ or $d = ec$. By symmetry we may assume $d=ec$. But  then the second row of $\hat M_{m-1}$ is a left multiple of the first row, which contradicts the minimality of the projective resolution.  This contradiction proves that any linearly dependent set of rows of $E_m$ must involve a row of zeroes.

Suppose now that some set of rows of $[L_m:E_m]$ is minimally linearly dependent.  By the above, the $E_m$ portion of every one of those rows must be zero.  We are left with a minimally linearly dependent set of rows of $L_m$.  But as above, the nonzero rows of $L_m$ are identical to the corresponding rows of $\hat M_m$, which has linearly independent rows.  Therefore at least one of the rows in our minimal dependence set is zero.  This proves the required claim. 
\end{pf}

The monomial algebra $A$ will fail to be $\K_2$ if and only if for some minimal $i$ there exists some $b\in S_i$ and $a\in L(b)$ such that $|a|>1$ and $ab\notin R$.    In this case we will say that the $\K_2$ property fails at level $i$.  

 \begin{ex}  Let $A$ be the monomial algebra generated by $\{w,x,y,z\}$ with relations $R=\{z^2y^2,y^3x^2,x^2w,zy^3x\}$.   Then $S_1=\{w, x,y,z\}$, $S_2=\{z^2y,y^3x, x^2, zy^3 \}$, $S_3=\{y^3 \}$, $S_4=\{z^2\}$ and $S_5=\emptyset$.   In this case  $\K_2$ does not fail at level 2, but $\K_2$ fails at level 3 because in $S_3$ we have $y^3$ which is minimally annihilated by $z^2$, however $z^2y^3$ is not in $R$.    
 \end{ex}
 
 \begin{ex}   Let $A$ be generated by $\{x,y,z\}$. Let $R=\{x^4, yx^3, x^3z\}$.   Then $S_1=\{x,y,z\}$, $S_2=\{x^3, yx^2 \}$ and $S_3 =\emptyset$.   We conclude that $A$ is $\K_2$
\end{ex}

For any monomial algebra if $b\in S_1$ and $a\in L(b)$ then $ab\in R$, so the earliest $\K_2$ can fail is at level two.  In the case of $N$-homogeneous algebras this is also the latest $\K_2$ can fail.

\begin{cor}
Let $A$ be an $N$-homogenous algebra.  $S_3$ will be empty if $A$ is $\K_2$,  and otherwise $\K_2$ will fail at level 2.  
\end{cor}  

\begin{pf}
Assume  that every monomial in $R$ is degree $N$.  Since every quadratic monomial algebra is Koszul, and hence $\K_2$, we will assume $N>2$.   Let $S_1=\{x_1,...,x_n\}$ and let $b$ be in $S_2$.  
Since $A$ in $N-$homogenous, $|b|=N-1$.   For any $a\in L(b)$ if $|a|>1$ then $|a|>N$, which means $ab\notin R$ and $\K_2$ fails at level 2.     If $|a|=1$ for all $a\in L(b)$  and all  $b\in S_2$, then $|ab|=N$ and $\overline{ab}=0$, which means $ab\in R$.  In this case $\K_2$ has not failed at level 2 and moreover, since $|a|=1$ we have $a\in S_1$ so that $S_3=\emptyset$.  If  $S_3=\emptyset$ then $S_m=\emptyset$ for all $m>2$     and $\K_2$ does not fail at any level.  
\end{pf}

This Corollary should be compared to the overlap condition for $N$-Koszul algebras given in \cite{Berg} Proposition 3.8.
The two conditions are logically equivalent, but the condition here requires fewer formal calculations.

%%%%%%%%%%%%%%%%%%%%%%%%%%%%%%%%%%%%%
\section{Spectral sequence lemmas}

This section  contains some spectral sequence facts which will be used in the subsequent sections.  Let $\phi:A\to B$ be a (graded) ring epimorphism with kernel $J$.  For any $B$ module $M$ we write 
$^\phi M$ for the corresponding $A$ module induced via $\phi$. For any graded $A$ module $N$ we write
$N^J$ for the graded $B$-module $N^J = \{n\in N\ |\ Jn=0\} = Hom_A(B,N)$.
The following lemma about the change of rings spectral sequence of Cartan and Eilenberg is certainly well known to experts.  

\begin{lemma}\label{module1}
Let $\phi:A \onto B$ be an epimorphism of $K$ algebras.    For any left $A$-module $N$ and any 
left $B$-modules $M$ and ${'M}$ we have first quadrant spectral sequences $E_r$ and 
$'E_r$, with $E_2$-terms  $E_2^{p,q} = Ext_B^p(M,Ext_A^q(B,N))$ and $'E_2^{p,q} = Ext_B^p({'M},Ext_A^q(B,N))$ which converge to $Ext_A({^\phi M},N)$ and $Ext_A({^\phi{'M}},N)$ respectively.  The differentials $d_r$ and $'d_r$ of these spectral sequences intertwine with the natural cup product by 
$Ext_B('M,M)$ from the right, that is we have the commutative diagram:
$$\begin{array}{ccc} E_r \times Ext_B('M,M) & \buildrel d_r \times 1\over \longrightarrow & E_r \times Ext_B('M,M)\\[6 pt]
\hphantom{\cupp}\downarrow \cupp && \hphantom{\cupp}\downarrow \cupp\\
'E_r & \buildrel 'd_r \over \longrightarrow & 'E_r\\[6 pt]
\end{array}
$$

\end{lemma}

\begin{pf}  We prove the lemma only for $r=2$, the general case being similar. We recall the construction of the spectral sequences.  Let $Q\to M\to 0$ and $'Q\to {'M}\to 0$ be projective resolutions (over $B$) of $M$ and $'M$ respectively.  Let $0\to N \to I$ be an $A$-injective resolution of $N$. We write $J_q = Hom_A(B,I_q)$. The spectral sequence $E_r$ is the first spectral sequence related to the double complex $F^{p,q} = Hom_B(Q^p,Hom_A(B,I_q)) = Hom_B(Q^p,J_q)$ with horizontal differential $d_h$ and vertical differential $d_v$ (chosen so that $d_h d_v + d_v d_h = 0$).  Similarly we have 
$'F^{p,q} = Hom_B('Q^p,J_q)$ with differentials $'d_h$ and $'d_v$.  We can use the usual trick (see for example \cite{Weibel}) of identifying $E_2^{p,q}$ with the set of classes of elements $(f,g)\in F^{p,q}\times F^{p+1,q-1}$ for which 
$d_vf=0=d_hf+d_vg$, modulo elements of the form $(0,g)$, $(d_vf,d_hf)$ 
(for $f\in F^{p,q-1}$) or $(d_hg,0)$ (for $g\in F^{p-1,q}$).  This allows us to identify $d_2:E_2 \to E_2$ as $d_2[(f,g)] = [(d_hg,0)]$.  

Let $[\zeta]\in Ext_B^k('M,M)$ be represented by $\zeta \in Hom_B('Q^k,M)$.  The following diagram then defines, via projectivity, maps $\zeta_j$, $j\ge 0$, so that the top row of boxes are commutative and the box in the bottom row is anticommutative.
$$\begin{array}{ccccccccccccccccc}
'Q^{k+p+2}&\kern -10pt\to&\kern -10pt 'Q^{k+p+1} &\kern -10pt \to &\kern -10pt 'Q^{k+p} &\kern -10pt \to \cdots  \to &\kern -10pt 'Q^{k+1} & \kern -10pt\to & \kern -10pt 'Q^{k} \\
\zeta_{p+2}\downarrow\hphantom{\zeta_{p+2}} &&
	\kern -10pt\zeta_{p+1}\downarrow\hphantom{\zeta_{p+1}}&& 
	\kern -10pt\zeta_p\downarrow\hphantom{\zeta_{p}} &
	\cdots &\zeta_1\downarrow\hphantom{\zeta_{1}} &&
	\kern -10pt\zeta_0\downarrow\hphantom{\zeta_{0}}&
	\kern -10pt\buildrel \zeta \over \searrow\\[7 pt] 
Q^{p+2} &\kern -10pt \to& \kern -10pt Q^{p+1} &\kern -20pt \to &\kern -10pt Q^{p} & \kern -10pt \to \cdots \to & \kern -10pt Q^{1} & \kern -10pt\to &\kern -10pt Q^0&\kern -10pt \to & M \\[7pt]
&&g\downarrow\hphantom{g}&&f\downarrow\hphantom{F}\\[7pt]
&& \kern -10ptJ_{q-1} &\kern -10pt \to & \kern -10pt J_q 
\end{array}
$$
From this diagram it is clear that $[(f,g)] \cupp [\zeta] = [(f\zeta_p,g\zeta_{p+1})]$.
So we may calculate $'d_2([(f,g)] \cupp [\zeta]) = [('d_h(g\zeta_{p+1}),0)]
=[(d_h(g)\zeta_{p+2},0)] = [(d_hg,0)]\cupp [\zeta] = d_2[(f,g)] \cupp [\zeta]$, as required.
\end{pf}

The following lemma should also be well-known.

\begin{lemma}\label{module1A}
Let $\phi:A \onto B$ be an epimorphism of $K$ algebras.    For any left $A$-modules
$N$ and $'N$ and any 
left $B$-module $M$ consider the first quadrant spectral sequences $E_r$ and 
$'E_r$, with $E_2$-terms  $E_2^{p,q} = Ext_B^p(M,Ext_A^q(B,N))$ and $'E_2^{p,q} = Ext_B^p(M,Ext_A^q(B,{'N}))$, which converge to $Ext_A({^\phi M},N)$ and 
$Ext_A({^\phi{M}},{'N})$ respectively.  Then, the left cup product action of $Ext_A(N,{'N})$ from
$Ext_A(B,N)$ to $Ext_A(B,{'N})$, commutes with the natural action of $B$.  Hence
the $Ext_A(N,{'N})$-action lifts to an action (still denoted $\cupp$) from 
$E_r$ to $'E_r$.  Moreover, the differentials $d_r$ and $'d_r$ of these spectral sequences intertwine with this action. That is, we have the commutative diagram:
$$\begin{array}{ccc} Ext_A(N,{'N})\times E_r & \buildrel 1\times d_r \over \longrightarrow 
	&  Ext_A(N,{'N})\times E_r\\[6 pt]
\hphantom{\cupp}\downarrow \cupp && \hphantom{\cupp}\downarrow \cupp\\
{'E}_r & \buildrel 'd_r \over \longrightarrow & {'E}_r\\[6 pt]
\end{array}
$$
\end{lemma}

\begin{pf}
The first claim is clear (since the action of $B$ is by precomposition and the action of
$Ext_A(N,{'N})$ is by post composition).   Let $Q\to M\to 0$  be a projective resolution (over $B$) of $M$.  Let $0\to N \to I$ and $0\to {'N} \to {'I}$ be $A$-injective resolutions of $N$ and ${'N}$.   Let $F^{p,q} = Hom_B(Q^p,Hom_A(B,I_q)) = Hom_B(Q^p,I_q^J)$ with horizontal differential $d_h$ and vertical differential $d_v$.  Similarly let 
$'F^{p,q} = Hom_B(Q^p,Hom_A(B,{'I_q}))=Hom_B(Q^p,{'I_q^J})$ with differentials $'d_h$ and $'d_v$.  As in the previous proof, we identify an element of $E_2^{p,q}$ with a certain class of elements $(f,g)\in F^{p,q}\times F^{p+1,q-1}$.  Let $[\zeta]\in Ext_A^k(N,{'N})$ be represented by $\zeta \in Hom_A(N,{'I_k})$.  The following diagram defines, via injectivity, maps $\zeta_j$, $j\ge 0$ so that the bottom row of boxes is commutative and
the top box is anticommutative.  Note that we abuse notation and consider $f$ and 
$g$ as maps into $I_q$ and $I_{q-1}$ respectively.  
$$\begin{array}{ccccccccccccccccc}
&&&&&& \kern -10pt Q_{p+1} &\kern -10pt \to &\kern -10pt Q^p \\[7pt]
&&&&&&\kern -10pt g\downarrow\hphantom{g}&&\kern -10pt f\downarrow\hphantom{F}\\[7pt]
0\to N&\kern -5pt\longrightarrow&\kern -5pt I_0 &\kern -10pt \to &\cdots & \to &I_{q-1}&\kern -10pt  \to &\kern -10pt I_q &\kern -10pt \to &\kern -5pt I_{q+1} \\[7pt]
&\kern -5pt\zeta \searrow&\kern -5pt\zeta_{0}\downarrow\hphantom{\zeta_{0}} &&&&
	\kern -10pt\zeta_{q-1}\downarrow\hphantom{\zeta_{q-1}}&&
	\kern -10pt\zeta_q\downarrow\hphantom{\zeta_{q}}&&
	\kern -10pt \downarrow \\[7 pt] 
&&\kern -5pt{'I_k}  &\kern -10pt \to& \cdots & \to &\kern -10pt {'I_{k+q-1}} &\kern -10pt \to & \kern -10pt {'I_{k+q}}&\kern -10pt \to &\kern -5pt 'I_{k+q+1} \\[7pt]
\end{array}
$$
From this diagram we see that $[\zeta]\cupp[(f,g)] = [(\zeta_q f,\zeta_{q-1}g)]$.  (We are once again abusing notation, using the fact that image of $\zeta_qf$ must be in $'I_{r+q}^J$ and similarly for $\zeta_{q-1}g$.)  Thus we have
$'d_2([\zeta]\cupp[(f,g)]) = [('d_h(\zeta_{q-1}g),0)] = [(\zeta_{q-1}d_h(g),0)]
= [\zeta]\cupp [(d_h(g),0)] = [\zeta]\cupp d_2[(f,g)].$  This proves the lemma for $r=2$, and the general proof is essentially the same. 
\end{pf}

For the remainder of this section  
let $A$ be a connected graded $K$-algebra and $g\in A$ a homogeneous element of degree $d$ that is normal and regular in $A$.  Set $B=A/gA$ and $\phi:A \onto B$ the associated epimorphism. For any graded $A$-module $N$, let $N^g$ be the $B$-module 
$\{n\in N | gn=0\}$. Then the Cartan-Eilenberg spectral sequence for 
$Ext_A({^\phi M},N)$,   $E_2^{p,q} = Ext_B^p(M,Ext_A^q(B,N))$, satisfies 
$E_2^{p,0} = Ext_B^p(M,N^g)$, $E_2^{p,1} = Ext_B^p(M,N/gN)$ and $E_2^{p,q}=0$ for $q>1$.  The edge homomorphism $E_2^{p,0}\to Ext_A^p({^\phi M},N)$ 
of this spectral sequence is the map induced by 
$\phi^*$ and we denote it by $\phi^*$ as well.  Since the spectral sequence has only two nonzero rows, there is an associated long exact sequence:
$$\begin{array}{l} \cdots \to Ext_A^n({^\phi M},N) \buildrel \gamma \over \to Ext_B^{n-1}(M,N^g) 
\buildrel d_2 \over \to Ext_B^{n+1}(M, N/gN)  \\[10 pt]
\hfill \buildrel \phi^* \over \to Ext_A^{n+1}({^\phi M},N) 
\buildrel \gamma \over \to Ext_B^n(M,N^g) \to \cdots\ \ 
\end{array} $$ 
We are primarily interested in applying this long exact sequence to the case when $M$ and $N$ are both the trivial module, $K$. But we need to know that the map $\gamma$, like $d_2$, intertwines with the right cup product from $Ext_B('M,M)$.  This theorem should also be well-known and we omit the proof. 

\begin{thm}\label{module2}
Let $A$, $g$, $B$ and $\phi$ be as above.   For any left $A$-module $N$ and any 
left $B$-modules $M$ and $'M$ let  $E_2^{p,q}$ and ${'E}_2^{p,q}$ be as in Lemma 
\ref{module1}.  There is a well defined cup product:
$$Ext_A^n(^\phi M,N) \times Ext_B^m('M,M)\buildrel \cupp\over \to 
Ext_A^{n+m}(^\phi {'M},N).$$
Moreover,  the maps 
$\gamma: Ext_A^{p+1}({^\phi M},N)\to E_2^{p,1}$ and ${'\gamma}:Ext_A^{p+1}({^\phi{'M}},N) \to {'E}_2^{p,1}$ intertwine with the right action of $Ext_B({'M},M)$.  That is we have a commutative diagram
$$\begin{array}{ccc} Ext_A^p({^\phi M},N) \times Ext_B^k('M,M) & \buildrel \gamma\times 1 \over \longrightarrow & E_2^{p,1} \times Ext_B^k('M,M)\\[6 pt]
\hphantom{\cupp}\downarrow \cupp && \hphantom{\cupp}\downarrow \cupp\\
Ext_A^{p+k}({^\phi{'M}},N) & \buildrel {'\gamma} \over \longrightarrow & {'E}_2^{p+k,1}\\[6 pt]
\end{array}
$$

\end{thm}

Consider now the case $N=M=K$, the trivial module for either $A$ or $B$ (from now on we will suppress the notation $^\phi K$).  The spectral sequence above then becomes
$E_2^{p,q}=0$ for $q>1$, $E_2^{p,1} = Ext_B^p(K,K(d)) = E^p(B)(d)$ and 
$E_2^{p,0} = Ext_B^p(K,K) = E^p(B)$.  The associated long exact sequence is therefore
\begin{equation}\label{spectral1}
  \cdots \to E^{p-2}(B)(d)\buildrel d_2\over \to E^p(B) \buildrel \phi^*\over \to
 E^p(A) \buildrel \gamma \over \to E^{p-1}(B)(d) \to \cdots 
 \end{equation}
 The lemmas show that this can be interpreted as an exact triangle of {\it right} $E(B)$ modules.  
 
 There is a second Cartan-Eilenberg spectral sequence converging to $E(A)$ given by
 $\bar E_2^{p,q} = Ext^p(Tor_q^A(B,K),K)$.  This spectral sequence is also supported on $p\ge 0$ and $q=0,1$ with nonzero terms $\bar E_2^{p,1} = Ext_B^p(K(-d),K) = E^p(B)(d)$ and $\bar E_2^{p,0} = Ext_B^p(K,K)= E^p(B)$.  The edge homomorphism of this spectral sequence is also $\phi^*: E(B) \to E(A)$ and we have a long exact sequence:
\begin{equation}\label{spectral2}
 \cdots \to E^{p-2}(B)(d)\buildrel \bar d_2\over \to E^p(B) \buildrel \phi^*\over \to
E^p(A) \buildrel \bar\gamma \over \to E^{p-1}(B)(d) \to \cdots
\end{equation}
Exactly as above, we see that the differential $\bar d_2$, as well as the induced map
$\bar \gamma$ intertwine with the {\it left} cup product by elements of $E(B)$.   We therefore interpret
\eqref{spectral2} as an exact triangle of {\it left} $E(B)$-modules.

Finally, we note that although the terms of the two exact triangles above are all the same, and the map $\phi^*$ is common to them,  the sequences are not, in general, the same.  In particular, they are not generally triangles of $E(B)$-bimodules. 
 
%%%%%%%%%%%%%%%%%%%%%%%%%%%%%%%%%%%%%%%%%%

\section{$\K_2$ modules and a change of rings theorem}

For Koszul algebras, Positselskii \cite{Posit}  establishes a powerful change of rings theorem, which states that if $A$ is a quadratic algebra and $B$ is a quadratic factor algebra of $A$ which admits a linear free resolution as an $A$-module, then the Koszul property for $B$ lifts to the Koszul property for $A$, (see also \cite{Polish}).   In this section we prove  the 
$\K_2$ analog of that theorem using the definition of a $\K_2$ module given below.

Let $A=T(V)/I$ be a graded algebra and $W$ a graded left $A$-module.  Let 
$D(A)$ be the subalgebra of $E(A)$ generated by $E^1(A)$ and $E^2(A)$.  
Let $Q^* \to W \to 0$ be a projective resolution of $W$ as a left $A$ module, where the maps
$Q^{n} \to Q^{n-1}$ are given by matrices $M_n$ with entries in $T(V)$.  
Let $L_n$ be the linear part of $M_n$, i.e. $M_n$ modulo $T(V)_{\ge 2}$.  For 
$n\ge 2$, let $E_n = M_nM_{n-1}$ modulo $I'$ and let $E_1=0$. 

\begin{defn}  We say $W$ is a $\K_2$ $A$-module if for all $1\le n \le pd_A(W)$, the 
rows of the matrix $[L_n:E_n]$ are linearly independent.  
\end{defn}

Theorem \ref{MC} can be interpreted as the statement that the algebra $A$ is $\K_2$ if and only if the trivial module is a $\K_2$-module.  The following lemma is proved in exactly the same way as that theorem.

\begin{lemma}\label{K_2-module}  The module $W$ is $\K_2$ if and only if $Ext_A(W,K)$ is generated as 
a left $D(A)$ module by $Ext_A^0(W,K)$.  
\end{lemma}

Before we can state and prove the analog of Positselskii's theorem for $\K_2$ algebras, we need a technical lemma. Let $A=T(V)/I$ be a graded algebra and let $B$ be a graded factor algebra of $A$. Let
$\phi:A\to B$ be the associated epimorphism.  

\begin{lemma}  Suppose that $B$, as a left $A$-module, is a $\K_2$-module.  Then the natural action of  $B$ on $Ext_A(B,K)$ is trivial.
\end{lemma}

\begin{pf}
Assume that  the 
rows of the matrix $[L_n:E_n]$ are linearly independent and fix $b\in B_k$ for some $k>0$.  The left action of $b$ on $Ext_A(B,K)$ is induced from the right action of $b$ on $B$.  We choose matrices $N_i$ with homogeneous entries in $T(V)$ to produce a commutative diagram of $A$-modules:
$$
\begin{array}{ccccccccccc}
\to& Q^i &\kern -16pt \buildrel M_i\over \to&  Q^{i-1}&\kern -16pt  \to& \cdots &\to &Q^0 &\kern -16pt\to& \kern -16pt B &\kern -16pt \to 0\\
& N_i\downarrow\hskip 20pt & & N_{i-1} \downarrow\hskip 30pt &&&& N_0\downarrow\hskip 20pt &&\kern -16pt \cdot b\downarrow\hskip 10pt \\
\to& Q^i &\kern -16pt \buildrel M_i\over \to&  Q^{i-1}& \kern -16pt \to& \cdots &\to &Q^0 &\kern -16pt\to& \kern -16pt B &\kern -16pt \to 0\\
\end{array}
$$
To prove the lemma, it suffices to prove the claim that there are no non-zero entries in $N_i$ of degree 0.  We do this by induction on $i$.  Write 
$N_i = N_i^0 + N_i^+$, where the entries of $N_i^0$ are of degree 0 and the entries of 
$N_i^+$ are of positive degree.  We must prove $N_i^0 = 0$ for all $i>0$. 

For $i=0$, $N_0$ is a diagonal matrix with $b$ on the diagonal, hence $N_0^0=0$.   The commutativity of the above diagram implies that $N_1M_1= M_1 b\ mod\ I$.  Since $I$ has no elements of degree 1 the linear part of $N_1M_1$ is the linear part of $M_1 b$, which is $0$.  $N_1^0L_1$ is the linear part of $N_1M_1$, hence $N_1^0L_1=0$. Since the rows of $[L_1:E_1]$ are linearly independent, $N_1^0=0$.  Assume now that $i>1$ and the claim is true for
all $k<i$.   Consider the matrix $N_i^0[L_i:E_i] = [N_i^0L_i:N_i^0E_i]$.  First we note that $N_i^0L_i$ is the linear component of the matrix $N_iM_i$.  Since $I$ contains no elements of degree 1  this is also the linear component of $M_iN_{i-1}$.  By the inductive hypothesis, every entry of this matrix has degree at least 2.  Thus $N_i^0L_i = 0$.  Similarly, $N_i^0E_i = N_i^0(M_iM_{i-1} mod\  I') = 
(N_i^0M_i M_{i-1} mod\  I')$.  By definition, $(N_i^+M_iM_{i-1} mod\  I') = 0$, since every entry of $M_iM_{i-1}$ is already in $I$.  
Now $N_iM_i=M_iN_{i-1}\ mod\ I$, so $N_iM_iM_{i-1}=M_iN_{i-1}M_{i-1}\ mod\ I'$, and likewise  $M_iN_{i-1}M_{i-1}=M_iM_{i-1}N_{i-2}\ mod\ I'$.
Thus, again by induction,  $N_i^0E_i = (N_iM_iM_{i-1}mod\  I') = (M_iM_{i-1}N_{i-2} mod\  I') = 0$.  Hence $N_i^0$ annihilates the matrix $[L_i:E_i]$.  Since the rows of $[L_i:E_i]$ are linearly independent it follows that
$N_i^0 = 0$.  This completes the induction and proves the lemma.
\end{pf}

\begin{thm}\label{ringchange}  Let $A$ be a graded algebra and $B$ a graded factor algebra of $A$.  Assume that $B$ is $\K_2$ as an $A$-module and also $\K_2$ as an algebra.  Then
$A$ is $\K_2$.  
\end{thm}  

\begin{pf}
Let $\phi: A\to B$ be the algebra epimorphism.  As usual, we have a change of rings spectral sequence $E_r^{p,q}$ converging to $E(A)$ with $E_2$ term
$E_2^{p,q} = Ext_B^p(K,Ext_A^q(B,K))$.  By Lemmas \ref{module1} and \ref{module1A} this is a spectral sequence of right $E(B)$ modules and left
$E(A)$ modules.  By the previous lemma, the $B$-action on $Ext_A(B,K)$ is trivial, and
thus we have $E_2^{p,q} = Ext_A^q(B,k)\otimes_K E(B)$ as an $E(A)$-$E(B)$ bimodule.
Since $_AB$ is a $\K_2$-module, this bimodule is generated by 
$E_2^{0,0}=Ext_B^0(B,K)\otimes E^0(B)$ as a $D(A)$-$E(B)$ bimodule. 
Since $d_2$ is a bimodule homomorphism and $d_2(E_2^{0,0})=0$, we conclude that $d_2=0$.  Similarly, for all
$r\ge 2$, $d_r=0$ and $E_r^{p,q} = E_2^{p,q}$.  

Since our spectral sequence has collapsed, there is a filtration $F^rE(A)$ on 
$E(A)$ for which $(F^{q}E^{p+q}(A)/F^{q-1}E^{p+q}(A)) = E_2^{p,q}$.  $E(A)$ is a
$D(A)$-$E(B)$-bimodule, where $E(B)$ acts on the right through $\phi^*$.  The convergence of the spectral sequence is compatible with the two bimodule actions.
Hence we conclude that $E(A)$ is generated by $E^0(A)$ as an $D(A)$-$E(B)$ bimodule.  By hypothesis $\phi^*(E(B)) = \phi^*(D(B)) \subset D(A)$, and thus
$E(A)$ is generated by $E^0(A)$ over $D(A)$, i.e. $E(A)=D(A)$, as required.   
\end{pf}
 
\begin{example}  Let $A$ and $B$ be defined by
$$A=\frac{K\langle x,y,z\rangle} {\langle xz-zx, yz-zy, x^3z, y^4+xz^3\rangle}, \ \ B=A/\langle z \rangle=\frac{K\langle x,y\rangle} {\langle  y^4\rangle}$$ 
 
Note that $A$ is not a regular central extension of $B$. Clearly $B$ is a $\K_2$ algebra, and we claim $B$ is also a $\K_2$ $A$-module.   To see this it suffices to show that the following is a projective resolution of $B$ as an $A$-module.  We omit any details.
$$ \begin{array}{r}\cdots \to A(-8) {\buildrel
	\left(\begin{array}{c}x^3\end{array}\right)
	\over \longrightarrow} A(-5){ \buildrel \left(\begin{array}{c}z \end{array}\right) \over \longrightarrow}
A(-4) {\buildrel
	\left(\begin{array}{c}x^3\end{array}\right)
	\over \longrightarrow}
A(-1)  {\buildrel\left(\begin{array}{c} z\end{array}\right)\over \longrightarrow} 
A \to B \to 0.
\end{array}$$
We conclude that $A$ is also a $\K_2$ algebra.  
\end{example}

%%%%%%%%%%%%%%%%%%%%%%%%%%%%%%%%%%%%%%%%%% 
\section{Normal, regular factor rings I, the degree 1 case}

Let $A=T(V)/I$ be a graded $K$-algebra and $g$ a homogeneous element of degree 1 that is both normal and regular in $A$.  We set
$B=A/gA$ and let $\phi:A\to B$ the natural graded algebra epimorphism.  The purpose of this section is to investigate the extent to which $A$ and $B$ inherit the $\K_2$ property from one another.  We remind the reader that the definition of $\K_2$ was specifically motivated by the desire that the class of algebras be closed when passing from 
$B$ to $A$.  Fortunately, this holds as a consequence of Theorem \ref{ringchange}.

\begin{thm}    If the algebra $B$ is $\K_2$ then the algebra $A$ is $\K_2$.   
\end{thm}
 
The converse of the previous theorem, as we will see, is not true.  However, one extra hypothesis will give us a converse.

\begin{thm}\label{quotient>1}  Let $A$ be $\K_2$ and assume that $E^3(B)$ is generated by $E^1(B)$ and $E^2(B)$.  Then $B$ is $\K_2$.  
\end{thm}

\begin{pf}
Assume $A$ is $\K_2$. Consider the exact triangle of right $E(B)$-modules $\eqref{spectral1}$, with $d=1$.
Since $E^0(B)(1)$ is supported in degree $-1$ and $E^{2,1}(B)=0$, the map
$d_2$ must be zero, the ring homomorphism $\phi^*:E(B)\to E(A)$ is injective and $\gamma$ is surjective. Choose any $\hat g\in E^{1,1}(A)$ for which $\gamma(\hat g) = 1\in E^0(B)$. Since the image of $\gamma$ is a projective (bigraded-free) right $E(B)$-module, we can split what is now a short exact sequence to get
 $E(A) = E(B) \oplus \hat g \cupp E(B)$  (we supress the $\phi^*$-notation).  Now consider the long exact sequence \eqref{spectral2}.  Since 
 $E^1(A) = E^1(B) \oplus K\hat g$ and $E^1(B)$ is in the kernel of $\gamma'$,  
 we see $\gamma'(\hat g) \ne 0$. Hence, arguing as above, we get $E(A) = E(B) \oplus E(B) \cupp \hat g$.   
 
 Let $D = \sum_{j,k}D^{j,k}$ be the bigraded subalgebra of  $E(B)$ generated by  $E^1(B)$ and $E^2(B)$.  
We have $E^1(A) = D^1 + K\cupp \hat g$ and $E^2(A) = D^2 \oplus D^1 \cupp \hat g = D^2 \oplus \hat g\cupp D^1$.  Thus $D$ and $\hat g$ generate $E(A)$.  Since $E^3(B) = D^3$, we have 
$D^2\cupp \hat g \subset D^3 \oplus \hat g\cupp D^2$.  Therefore $D\cupp \hat g \subset D\oplus \hat g \cupp D$.  Since $D$ and $\hat g$ generate $E(A)$, we conclude that 
$E(A) = D\oplus \hat g\cupp D$.  This proves that $D = E(B)$, as required.
\end{pf}

\begin{rmk}  The hypothesis that $E^3(B)$ should be generated by $E^1(B)$ and 
$E^2(B)$ is easily translated into the hypothesis that the matrix $[L_3:E_3]$, from Theorem \ref{MC}, has linearly independent rows. 
\end{rmk}

The following example shows that the $E^3(B)$ hypothesis of the theorem cannot be avoided.  The example also highlights an interesting connection between Theorem \ref{quotient>1} and the main theorem of \cite{PBW}, where it is shown that the regularity of a central extending variable is  controlled by information encoded in $E^3(B)$.  

\begin{ex}   Let $A=K\langle x,y,z,w,g\rangle/I$ where 
$I$ is generated by 
$$\{y^2z, zx^2 + gw^2, y^2 w^2,xg-gx,yg-gy,wg-gw,zg-gz\}.$$ 
Let $B = A/gA$. We make the following three non-obvious claims:

(1)  $B$ is not $\K_2$.

(2)  The element $g$ is central and regular in $A$.  

(3)  $A$ is $\K_2$
\end{ex}

\begin{pf}  Let $T=K\langle x,y,z,w,\rangle$.  

To prove the first claim, we apply the algorithm of Section \ref{monom} to get a minimal projective resolution of the trivial $B$-module.  This has the form
$$ 0 \to B(-5) \buildrel M_3\over \to B(-3,-3,-4)\buildrel M_3\over \to B(-1,-1,-1,-1)\buildrel M_1\over \to B\to K\to 0$$
where 
$$M_3 = \left(\begin{array}{ccc}0 & y^2 & 0\end{array}\right),\quad
M_2 = \left(\begin{array}{cccc} 0&0&y^2&0\\ zx&0&0&0\\ 0&0&0&y^2w\end{array}\right),
\quad
M_1=\left(\begin{array}{c}x\\ y\\ z\\ w\end{array}\right)$$
We calculate $\left[ L_3:E_3\right]$ for this, as in theorem \ref{MC}, and see that its first row is zero.  This proves (1).

To prove (2) we use the algorithm from \cite{PBW} for checking the regularity of a central extending variable.  We work over the algebra $T[g]$.  Define:
$$f_2 = \left(\begin{array}{c}0\\ w^2\\ 0\end{array}\right)\ 
\hbox{ and }
f_3= \left(\begin{array}{cccc}0&0&0&0\end{array}\right)$$
so that $M_2M_1 + f_2g$ and $M_3M_2 - f_3g$ are both identically $0$ in $A$.
Put 
$$\hat M_n = \left(\begin{array}{cc} M_n & f_n\\ (-1)^{(n-1)} g I & M_{n-1}\end{array}\right)$$
for all $n\ge 1$.    The principal theorem of \cite{PBW} states that $g$ is regular in $A$ if and only if the matrix $\hat M_3 \hat M_2$ is zero in $A$.  Furthermore, if $g$ is regular in $A$, then the matrices $\hat M_n$ provide the maps in a minimal projective resolution of the trivial $A$-module.  We have
$$\hat M_3 =
\left(\begin{array}{ccccccc}0&y^2&0&0&0&0&0\\ g&0&0&0&0&y^2&0\\
0&g&0&zx & 0&0&0\\ 0&0&g&0&0&0& y^2 w\end{array}\right),
\hat M_2 = \left(\begin{array}{ccccc}0&0&y^2&0&0\\ zx&0&0&0&w^2\\ 0&0&0&y^2w&0\\
-g&0&0&0&x\\ 0&-g&0&0&y\\ 0&0&-g&0&z\\ 0&0&0&-g& w\end{array}\right).$$
It follows immediately that $\hat M_3\hat M_2$ is 0 in $A$ and that $g$ is regular. 

To prove the last claim, we consider the matrix condition \ref{MC}.  We will refrain from including the calculation.  The salient point is that the term $y^2w^2$, which is an essential relation of $A$, appears in the first row of $[L_3:E_3]$ and that 
$-g$ appears in the first row of $[L_4:E_4]$.  
\end{pf}

\begin{rmk}  The example above becomes substantially more complex if we make the rings $A$ and $B$ commutative (that is, factor $A$ and $B$ by the commutativity relations).  Nonetheless, using Theorem \ref{ringchange} it can be shown that the interesting features of the example remain the same, i.e. $A$ is $\K_2$, $g$ is regular, and $B$ is not $\K_2$.
\end{rmk}

\section{Normal, regular factor rings II, the degree $d>1$ case}

Throughout this section, let $A=T(V)/I$ be a graded algebra, $g$ a normal and regular homogeneous element of $A$ of degree $d>1$ and $B = A/gA$. As before, 
$\phi:A \to B$ is the associated algebra epimorphism.  We  examine the conditions under which the $\K_2$ property will descend from $A$ to $B$.  

Since the map $d_2:E(B)(d) \to E(B)$ associated to the long exact sequence 
\eqref{spectral1} is a right $E(B)$ module homomorphism, there is a
distinguished element $\hat g \in E^{2,d}(B)$ such that $d_2(a) = \hat g \cupp a$. 

\begin{thm}\label{normalreg}
Assume $A$ is $\K_2$ and the map $\gamma:Ext_A^2(K,K) \to Ext_B^1(K,K(d))$ of  \eqref{spectral1} is $0$.  Then $B$ is $\K_2$, $\hat g$ is a normal and regular element of 
$E(B)$ and $E(A) = E(B)/\hat g E(B)$.
\end{thm}

\begin{pf}

Since $d>1$, $\phi^*:E^1(B) \to E^1(A)$ is an isomorphism (i.e. $A$ and $B$ have the same set of generators). From \eqref{spectral1} the hypothesis on $\gamma$ tells us that
$\phi^*:E^2(B) \to E^2(A)$ is surjective.  

Let $D=\oplus_n D^n$ be the graded subalgebra of $E(B)$ generated by $E^1(B)$ and $E^2(B)$.  

We use induction on $n$ to prove the following two claims: $D^n=E^n(B)$ and 
$\phi^*:D^n \to E^n(A)$ is surjective.  We may assume the claim for $n-1$ and $n-2$.  Our long exact sequence \eqref{spectral1} then looks like:
$$0\to D^{n-2}(d)\buildrel d_2\over \to E^n(B)
\buildrel \phi^*\over \to E^n(A) \buildrel\gamma\over \to D^{n-1} \to \cdots$$
Since $A$ is $\K_2$, we have by induction, 
$$\begin{array}{rcl}E^n(A)& = &E^1(A)E^{n-1}(A) + E^2(A)E^{n-2}(A)\\
& =&
\phi^*(D^1)\phi^*(D^{n-1}) + \phi^*(D^2)\phi^*(D^{n-2})\\
& = &\phi^*(D^n)
\end{array}$$
This proves the second claim and also shows that $\gamma(E^n(A))=0$.  Since $d_2$ is multiplication by $\hat g \in E^2(B)=D^2$, we get $E^n(B) = D^n + \hat g\cdot D^{n-2} = D^n$, completing the induction.

We now see that $\phi^*$ is surjective, $\gamma = 0$ and $d_2$ is injective.  In particular, $\hat g$ is left regular in $E(B)$.  Since $\phi^*$ is surjective we also see, from 
\eqref{spectral2} that $\bar d_2$ is injective and $\bar \gamma=0$.   There exists 
$\hat g'\in Ext^{2,d}(B)$ for which $\bar d_2(a) = a\cupp \hat g'$.  Then
$\hat g \cupp E(B) = E(B)\cupp \hat g' = ker(\phi^*)$.  It follows that $\hat g = \mu\hat g'$ for some scalar $\mu\ne 0$ and hence $\hat g$ is normal and regular in $E(B)$, as required.
\end{pf}

There are several ways to assure the hypothesis that the image of 
$\gamma:Ext_A^2(K,K) \to Ext_B^1(K,K(d))$ is 0.  In particular, the hypothesis can only fail if there are relations of $A$ of degree $d+1$ (since $Ext_B^1(K,K(d))$ is supported in degree 
$d+1$.)  More specifically, let $\sigma$ be the automorphism of $A$ defined by 
$ag=g\sigma(a)$.  Then $\sigma$ defines an automorphism of $T(V)$ as well.  Choose 
$g' \in T(V)_d$, a preimage of $g\in A$. It is not difficult to see that the map fails to be zero if and only if there exists an {\it essential} relation of the form 
$x\otimes g' - g'\otimes \sigma(x)$ for some $x\in T(V)_1$.

\begin{cor}\label{complete}  If  $A$ is a commutative or graded-commutative $\K_2$-algebra and $g\in A_d$ is regular, then $A/gA$ is $\K_2$.   In particular, any graded complete intersection is $\K_2$. 
\end{cor}

The hypothesis $\gamma(Ext_A^2(K,K))=0$ is certainly not necessary in Theorem \ref{normalreg}.  For example, in the $\K_2$-algebra $A=K\langle x,y\rangle /\langle xy^2-y^2x\rangle$, $y^2$ is normal and regular and $A/y^2A = K\langle x,y\rangle/\langle y^2 \rangle$ is also $\K_2$, but $\gamma(Ext_A^2(K,K)) \ne 0$.  However, the theorem is false without the $\gamma$ hypothesis  as shown by the following important example.

\begin{ex}  Let $A = K\langle x,y\rangle/\langle x^2y - y x^2, xy^3 - y^3 x\rangle$ and let
$B = A/y^3A = K\langle x,y\rangle /\langle x^2y - yx^2, y^3\rangle$.  It is evident that $y^3$ is central in $A$, and by a straightforward application of Bergman's diamond lemma \cite{bergman} we see that  $y^3$ is  regular in $A$.  Moreover $A$ is $\K_2$.  This is proved by applying theorem \ref{MC} to the following minimal projective resolution of $_AK$:
$$ \begin{array}{r} 0 \to A(-5){ \buildrel \left(\begin{array}{cc} y^2 & -x\end{array}\right) \over \longrightarrow}
A(-3,-4) {\buildrel
	\left(\begin{array}{cc}yx & -x^2\\ y^3 & -xy^2\end{array}\right)
	\over \longrightarrow}
A(-1,-1)\\ {\buildrel\left(\begin{array}{c} x\\ y\end{array}\right)\over \longrightarrow} 
A \to K \to 0.
\end{array}$$
However, $B$ is not $\K_2$.  The minimal projective resolution of $B$ is
$$\begin{array}{r} \cdots \to B(-5,-4)
{\buildrel\left(\begin{array}{cc} y^2 & x^2\\ 0 & y\end{array}\right)\over \longrightarrow}
B(-3,-3) 
{\buildrel\left(\begin{array}{cc} y x & -x^2\\ 0 & y^2\end{array}\right)\over \longrightarrow}
B(-1,-1)\\
{\buildrel\left(\begin{array}{c} x\\  y\end{array}\right)\over \longrightarrow}
B \to K \to 0,
\end{array}$$
which fails the criteria of $\ref{MC}$.  Alternatively, simply note that $E^1(B)$ and 
$E^2(B)$ are supported in internal degrees $-1$ and $-3$ respectively, and thus they can not possibly generate the nonzero cohomology class in $E^{3,5}(B)$. 
\end{ex}

%%%%%%%%%%%%%%%%%%%%%%%%%%%%%%%%%%%%%%%%%%%
\section{Graded Ore Extensions}

Fix a graded algebra $B$, a graded automorphism $\sigma$ of $B$ and
a degree +1 graded $\sigma$-derivation $\delta$ of $B$ (i.e. 
$\delta(ab) = \delta(a)\sigma(b) +a \delta(b)$).  We let $A$ be the associated Ore extension $B[z;\sigma,\delta]$.  The extending variable $z$ is assumed to have degree 1 in $A$.   In this section we investigate when the $\K_2$ property can pass between $A$ and $B$. 

 Let $J$ be the two sided ideal $J=AB_+ = B_+A \subset A$ and put $C=A/J$, with factor morphism $\phi:A\to C$.  We note that $C$ is just the polynomial ring $C=K[\bar z]$, $\bar z = \phi(z)$.   Let $\alpha:B\to A$ be the inclusion homomorphism. 

Let $\zeta:A\to A$ be the left $A$-module homomorphism given by {\it right} multiplication by $z$. Let $\zeta$ also denote the induced map  $\zeta:C \to C$ (still  given by multiplication by $z$).  We have a short exact sequence of left $A$-modules:
$$0\to C(-1) \buildrel \zeta \over \to C \to K \to 0$$
This induces a long exact sequence 
$$ \cdots\to Ext_A^{q-1}(C,K)(1) \to E^q(A) \to Ext_A^q(C,K) \buildrel \zeta^*\over \to Ext_A^q(C,K)(1) \to \cdots$$
Since $A$ is free as a right (or left) $B$-module and 
$_AC = A\otimes_B K$, $Ext_A(C,K) = Ext_B(K,K) = E(B)$, in particular we may consider
$\zeta^*$ as a map on $E(B)$.   As usual, let $\sigma^*$ be the automorphism of $E(B)$ induced from the automorphism $\sigma$. 

The following  Lemma does not seem to be widely known.

\begin{lemma}\label{cool}  The map $\zeta^*:E(B)\to E(B)$ is a $\sigma^*$-derivation (with respect to the cup product $\cupp$) of homological degree 0 and internal degree +1.  Moreover, $\zeta^*$ vanishes on  $E^1(B)$ and $E^2(B)$.
\end{lemma}

\begin{pf}  Let $Bar(B) = \bigoplus_n B\otimes B_+^{\otimes n}$ be the usual bar resolution of $_BK$, with differential 
$$\partial(b_0\otimes \cdots\otimes b_n) = \sum\limits_{i=1}^n (-1)^i b_0\otimes \cdots b_{i-1}b_i \otimes \cdots \otimes b_n.$$
Since $A_B$ is free and $A\otimes_B K = C$, we may tensor $Bar(B)$ by $A$ to get 
an $A$-projective resolution $A\otimes_B Bar(B) = \bigoplus_nA\otimes B_+^{\otimes n}$
of $_AC$ with boundary map $\partial_A$.
  
We define maps $\check\delta_n$ and $\zeta_n$ as follows:
$$\check\delta_n:B_+^{\otimes n} \to B_+^{\otimes n} \qquad \check\delta_n = \sum\limits_{i+j=n-1}1^{\otimes i}\otimes \delta \otimes \sigma^{\otimes j},$$
and
$$\zeta_n:A\otimes B_+^{\otimes n} \to A\otimes B_+^{\otimes n} \qquad \zeta_n = \zeta\otimes \sigma^n + 1\otimes \check\delta_n,\quad \zeta_0=\zeta.$$

It is straightforward to check 
$\zeta_n\circ\partial_A = \partial_A\circ \zeta_{n+1}$.   Therefore $\zeta^*$ on 
$Ext_A^m(C,K)$ is given by $[f] \mapsto [\zeta_q^*f]$, for  $f\in Hom_A(A\otimes B_+^{\otimes m},K)$ (a cocycle).  However, $(\zeta\otimes \sigma^{\otimes m})^*f = 0$, simply because $z$ has positive degree, and thus $\zeta_m^*f = (1\otimes \check \delta_m)^*f$.  In particular, under the identification
$Hom_A(A\otimes B_+^{\otimes m},K) = Hom_K(B_+^{\otimes m})$, $\zeta_m^*$ becomes 
$\check \delta_m^*$.

The definition of $\check \delta_{n+m}$ makes it clear that for $f\in Hom(B_+^{\otimes n},K)$ and 
$g\in Hom(B_+^{\otimes m},K)$, 
$\check\delta_{n+m}^*(f\otimes g) = f\otimes \check\delta_m^*(g) + \check\delta_n^*(f)\otimes (\sigma^{\otimes m})^*(g)$.   This proves that $\zeta^*$ is a $\sigma^*$-derivation, as claimed. 

Since $\zeta^*$ is a $\sigma^*$-derivation of internal degree $+1$, It is now clear that $\zeta^*$ vanishes on $E^1(B) = E^{1,1}(B)$ and on $E^{2,2}(B)$ (since $E^{2,1}(B)=0$).  Suppose now that 
$[f] \in E^{2,q}(B)$ for some $q>2$ and consider 
$\zeta^*[f] \in E^{2,q-1}(B)$.    Define $g\in Hom(B_+,K)$ as follows: $g(B_1)=0$ and for any  $b_i,c_i\in B_+$, $g( \sum_i b_ic_i) = f(\sum_i(\delta(b_i) \otimes \sigma(c_i) + b_i \otimes \delta(c_i))$.  From the formula $\partial(\sum_i(\delta(b_i) \otimes \sigma(c_i) + b_i \otimes \delta(c_i)) = 
\delta(\sum_i b_ic_i)$ and the fact that $f$ is a cocycle, it follows that $g$ is well-defined. 
But then $\partial^*g = \check\delta_2^*f$, showing that $\zeta^*[f] = 0$, as required.
\end{pf}

\begin{thm}\label{Ore}  If $B$ is a $\K_2$ algebra, then the graded Ore extension $A=B[z;\sigma,\delta]$ is also 
$\K_2$. 
\end{thm}

\begin{pf}

By the lemma, since $B$ is $\K_2$, the map $\zeta^*:E(B) \to E(B)$ is zero.  The long exact sequence
above becomes the short exact sequences
$$0 \to E^{q-1}(B)(1) \to E^q(A) \to E^q(B) \to 0$$
It is easy to see that the second map in this sequence is just the ring homomorphism $\alpha^*$.  Let $\check z\in E^{1,1}(A)$ span the kernel of $\alpha^*:E^1(A) \to E^1(B)$.  Then the fact that
$\alpha^*$ is surjective and has a kernel isomorphic to $E(B)(1)$ implies that $\check z$ generates the kernel of $\alpha^*$, i.e. our short exact sequence is encoding:
$$0\to \check z\cupp E(A) \to E(A) \to E(B) \to 0.$$
It follows at once that $A$ is $\K_2$.  
\end{pf}

\begin{rmk} If we hypothesize that the map $\zeta^*$ of Lemma \ref{cool} is known to be zero, then 
$\alpha^*:E(A)\to E(B)$ is surjective.  It follows then that $B$ inherits the $\K_2$ property from $A$.  
However, the following example shows that the map $\zeta^*$ need not be zero.  In this example, neither the algebra $B$ nor its Ore extension $B[z;\delta]$ are $\K_2$.   We do not know if the converse to \ref{Ore} is true or false.
\end{rmk}

\begin{example}  Let $B$ be the monomial algebra $K\langle x,y \rangle/\langle xyx, xy^2x,y^3\rangle$.  Let $\sigma$ be the identity automorphism of $B$ and let $\delta$ be the derivation of $B$ extended from the formula $\delta(x)=0$ and $\delta(y) = y^2$.    The algebra $B$   fails the 
$\K_2$ criterion of Section \ref{monom}.   
Since $B$ is a monomial algebra, each graded component $B_n$ has a canonical monomial basis, and these determine monomial bases for the graded components of $B_+^{\otimes m}$.  We define
$f\in Hom(B_+^{\otimes 3},K)_{-6}$ by: 
$$f(xy^2\otimes xy\otimes x) = f(xy\otimes yxy \otimes x) = f(x\otimes y^2xy \otimes x)=f(xy^2\otimes x\otimes yx) = 1$$ 
and by insisting that $f$ vanishes on all other monomials of degree 6.  One checks easily that
$\partial^* f = 0$, so that $[f]$ is a nonzero cohomology class in $E^{3,6}(B)$.  By the lemma,
$\zeta^*[f] = [\check \delta_3^*f] \in E^{3,5}(B)$.   We calculate:
$\check \delta_3^*f(xy\otimes xy\otimes x) = f(xy^2\otimes xy \otimes x + xy \otimes xy^2 \otimes x) = 1$
and similarly 
$\check \delta_3^*f(x\otimes yxy\otimes x) = f(x\otimes y^2xy \otimes x + x \otimes yxy^2 \otimes x) = 1$
and
$\check \delta_3^*f(xy\otimes x\otimes yx) = f(xy^2\otimes x \otimes yx + xy \otimes x \otimes y^2x) = 1.$
We see that $\check \delta_3^*f$ vanishes on all other monomials of degree 5.  This assures that
$[\check \delta_3^*f]\ne 0$ and shows that $\zeta^*:E^3(B) \to E^3(B)$ is not zero.  

\end{example}

\bibliographystyle{amsplain}
\bibliography{bibliog}

\end{document}